\documentclass[10pt,twocolumn]{article}
\usepackage{Nsc2016}
\usepackage{multirow}
\title{\uppercase{SIMULATION OF RECURSIVE FUNCTIONS BY MEANS OF INTERVAL ANALYSIS AND PSEUDO-ORBITS}}

\author[1]{\underline{Rodrigues Junior, H. M.}}
\author[2]{Peixoto, M. L. C.}
\author[3]{Nepomuceno, E. G.}

\affil[1]{Modelling and Control Research Group - GCOM, Federal University of S\~ao Jo\~ao del-Rei - UFSJ, \\  
    
    S\~ao Jo\~ao del-Rei - MG, Brazil heitormrjunior@gmail.com}
\affil[2]{marciapeixoto93@hotmail.com}
\affil[3]{nepomuceno@ufsj.edu.br}

\begin{document}

\maketitle

\Abstract{In this work we use intersection of different pseudo-orbits obtained by interval extensions to reduce the bounds of the exact solution provided by the toolbox Intlab. The method is applied on the logistic map.}

\Keywords{Modeling, Numerical Simulation and Optimization; Chaos and Global Nonlinear Dynamics; Bifurcation Analysis and Applications; Interval Analysis; Multiple Pseudo-Orbits.}

\section{INTRODUCTION}

Numerical computational techniques are almost never exact. This fact occurs due to limitations in representing real numbers by computers. Instead of a true result, its approximation is returned \cite{Galias2013}. Different software, machines and levels of precision may lead us to different results. Even for some mathematical equivalent functions, small rounding errors can be propagated, and after some iterates the result is no longer significant \cite{Galias2013,Nep2014,Nepomuceno2016}.

The error propagation problem is especially important for recursive functions that characterize chaotic systems, since small errors introduced in each computational step grow due the fact that the system is not contractive mapping \cite{Galias2013,HYG1987,Nep2014}. According to \cite{Nep2014}, a recursive function can be defined as
\begin{equation}
	\label{rec}
	x_{n}=f(x_{n-1}).
\end{equation}

A discrete-time series can be generated from a map by choosing some values for its parameters, an initial condition and iterating it recursively. The sequence of values of a map, represented by $ \{x_{0},x_{1},\dots,x_{n}\} $ is an orbit.

However, when a computer is used to iterate a map $ f $, we have a pseudo-orbit \cite{HYG1987}. In other words, there is an approximation of a true orbit due to inherent properties of digital computers and  we represent a specific pseudo-orbit $ i\in\mathbb{N} $ as $ \{\hat{x}_{i,0}, \hat{x}_{i,1},\dots,\hat{x}_{i,n}\} $, such that
\begin{equation}
\label{delta}
|x_{n}-\hat{x}_{i,n}|\leq\delta_{i,n},
\end{equation}
where $ \delta_{i,n}\geq0 $ and $ \delta_{i,n}\in\mathbb{R} $.

The logistic map \cite{May1976} is a recursive function given by
\begin{equation}
\label{log}
x_{n+1}=rx_{n}(1-x_{n}).
\end{equation}

This function describes a biological model in which $ x_{n} $ is a number between zero and one that represents the ratio of existing population to the maximum possible and $ r $ is the growth rate of the population. The values of interest for the parameter $ r $ are those between zero and four. May \cite{May1976} was the first to describe a rich set of dynamical properties subject to the variation to $ r $, varying from a fixed point to a chaotic behaviour.

\begin{figure}[h]
	\centering
	\includegraphics[width=9cm]{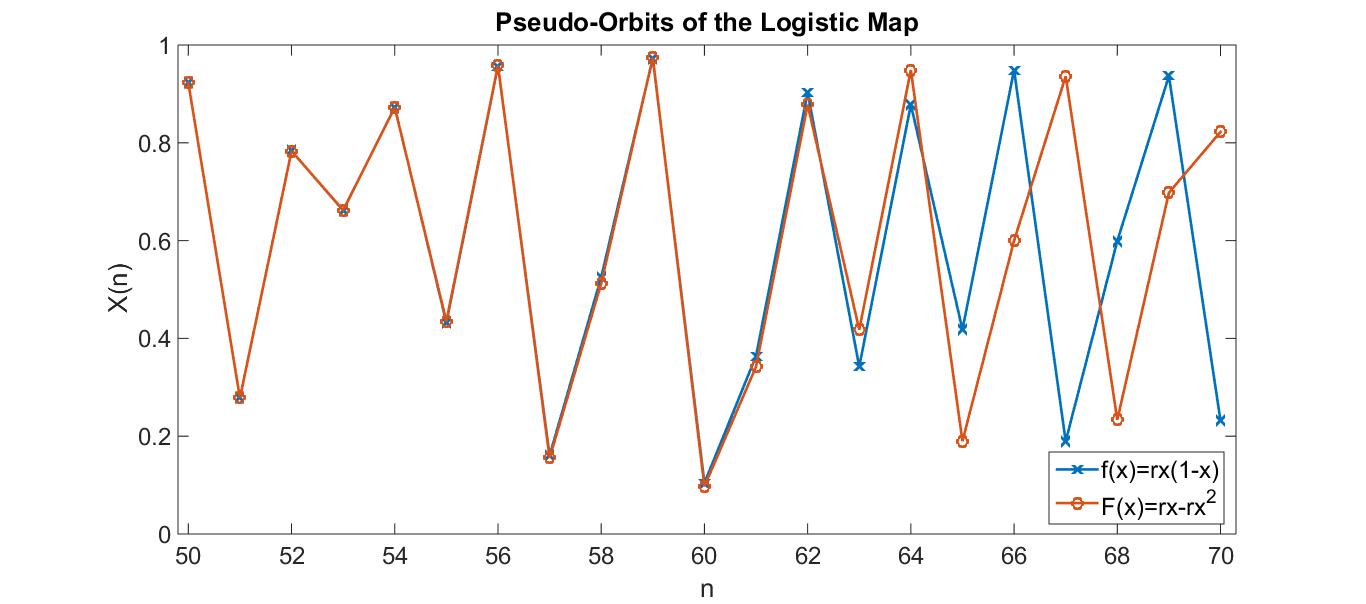}
	\caption{Simulation of (\ref{log}), with $ r=3.9 $ and $ x_{0}=0.6 $ for two different natural extensions. The simulation was performed in Matlab, which follows the IEEE 754-2008 standard for float point.}
	\label{sim1}
\end{figure}

To illustrate the error propagation, let us compute two pseudo-orbits, the first is given by $ f(x)=rx(1-x) $, and the second by its natural extension $ F(x)=rx-rx^{2} $, with $ r=3.9 $ and $ x_{0}=0.6 $. Any mathematical equivalent function of $ f $ is a natural extension. . The results obtained for the iterations 50 to 70 is shown in Figure \ref{sim1}. 

Clearly, the trajectory of these two pseudo-orbits is completely different after 60 iterations. For example, the distance between these pseudo-orbits at $ n=67 $ is greater than $ 0.74 $. Thus, both pseudo-orbits can not be accurate, and this lead us to a conclusion that even for very simple systems, propagation of rounding errors may lead to completely wrong results.

One way to try to reduce this problem is by means interval analysis \cite{MM1979}. While floating-point arithmetic is affected by rounding errors and can produce inaccurate results, interval arithmetic has the advantage of giving rigorous bounds for the exact solution. If the lower and upper bounds of the interval can be rounded down and rounded up, respectively, then finite precision calculations can be performed using intervals, to give an enclosure of the exact solution \cite{Ruetsch2005}.

\section{PURPOSE}

In interval analysis, real numbers are replaced by intervals. An interval is a closed set of real numbers $ x\in\mathbb{R} $ often denoted by an upper case letter $ X $ with endpoints $ \underline{x}\leq\overline{x} $ such that $ X=[\underline{x},\overline{x}]=\{x:\underline{x}\leq x\leq\overline{x}\} $.

Let $ f $ be a function of a real variable $ x $. An interval extension of $ f $ is an interval valued function $ F $ of an interval variable $ X $ with property
\begin{equation}
F(x)=f(x), \qquad x\in\mathbb{R}.
\end{equation}

Although mathematically equivalent, the interval extensions may have different sequences of arithmetic operations. Therefore it is important to mention that simple mathematical properties are not assured in interval arithmetic \cite{Nepomuceno2016}.

For a given interval $ X=[\underline{x},\overline{x}] $, its width is defined by $ w(X)=\overline{x}-\underline{x} $ and its centre is $ m(X)=0.5(\underline{x}+\overline{x}) $ \cite{alefeld2012}. The intersection of two intervals $ X\cap Y $ is a set of real numbers which belongs to both. The union $ X\cup Y $ is a set of real numbers which belongs to $ X $ or $ Y $ (or both). If $ X\cap Y $ is not empty, then 
\begin{equation}
\label{int}
w(X\cap Y)\leq \min\{w(X),w(Y)\}.
\end{equation}

The accuracy of a numerical computation in an interval is inversely proportional to the size of the interval. Therefore, the smaller the interval, the more accurate the information. Besides that, if $ X\cap Y\neq\emptyset $, the exact solution may belong to both intervals \cite{MM1979}.

The aim of this paper is to use the Intlab toolbox to show that the intersection of different pseudo-orbits obtained from interval extensions of the logistic map, generate a smaller interval where the real orbit is contained, increasing the accuracy of solution. 

The Intlab is a toolbox for Matlab supporting real and complex intervals designed to be very fast.  The Intlab concept splits in a fast interval library and an interactive programming environment for easy use of interval operations \cite{Rump1999}. 

\section{METHODS}

Let $ \{x_{0},x_{1},\dots,x_{n}\} $ be a true orbit and let $ [\hat{X}_{i,0},\hat{X}_{i,1},\dots,\hat{X}_{i,n}] $ be a pseudo-orbit obtained by some interval extension $ i\in\mathbb{N} $ such that  
\begin{equation}
\label{intd}
\hat{X}_{i,n}=[\underline{\hat{x}}_{i,n},\overline{\hat{x}}_{i,n}]=[x_{n}-\delta_{i,n},x_{n}+\delta_{i,n}].
\end{equation}

From (\ref{delta}) and (\ref{intd}) it is clear that
\begin{equation}
\label{pseudo}
x_{n}\in \hat{X}_{i,n}.
\end{equation}

From this we may stablish the following for two interval extensions $ F_{1} $ and $ F_{2} $ of the same function.
\begin{lemma}
$ x_{n}\in[\hat{X}_{1,n}\cap \hat{X}_{2,n}], \quad n\in\mathbb{N} $.\\
\end{lemma}
$Proof. $\quad Let us assume $ x_{n}\notin[\hat{X}_{1,n}\cap \hat{X}_{2,n}], \quad n\in\mathbb{N} $. Hence $ x_{n}\notin\hat{X}_{1,n} $ or $ x_{n}\notin\hat{X}_{2,n} $ or both. According to (\ref{pseudo}), it means that $ \hat{X}_{1,n} $ or $ \hat{X}_{2,n} $ is not a pseudo-orbit, which is a contradiction. That completes the proof. \hfill $\square$ \\

Besides that, according to  (\ref{int}) the intersection of these interval extensions can be used to limit the bounds of the true orbit $ x_{n} $. To illustrate this fact we will use different interval extensions to calculate two fixed points of distinct periods of the logistic map. If $ f(x_{n})=x_{n-p} $ than $ x_{n} $ is a fixed point of period $ p $. Based on \cite{Rodrigues-Junior2015}, in interval arithmetic we have the following.
\begin{lemma}
If $X_{i,n}\cap X_{i,n-p}\neq\emptyset$, then $x^{*}\in X_{i,n}$ is a fixed point of period $ p $.\\
\end{lemma}	
$ Proof. $\quad Let us assume that $ x^{*} $ is a fixed point of period $ p $. It means that $ x^{*}\in X_{i,n} $ and $ x^{*}\in X_{i,n-p} $, $ i,n,p\in\mathbb{R} $. If $X_{i,n}\cap X_{i,n-p}=\emptyset$, it means that $ x^{*}\notin X_{i,n} $ or $ x^{*}\notin X_{i,n-p} $ or both, which is a contradiction. \hfill $\square$ \\

In this work, first we will use the traditional approach to calculate the fixed points by an interval extension $ F_{1} $ of (\ref{log}). Than we will compare these calculations with a method which the calculation is made by the intersection at each iteration of two interval extensions $ F_{1} $ and $ F_{2} $ of (\ref{log}). We will use $ F_{1}(X)=rX(1-X) $ and $ F_{2}(X)=r(X(1-X)) $.



These methods will be compared by the width of the intervals and the number of iterations to get to some fixed point. Based on \cite{Fei1978}, three cases of fixed points of (\ref{log}) will be analized:

\begin{itemize}
	\item Case 1: $ r=3.3 $ and $ x_{0}=0.6 $, which orbit may converge to fixed point of period 2.
	\item Case 2: $ r=3.47 $ and $ x_{0}=0.6 $, which orbit may converge to fixed point of period 4.
	\item Case 3: $ r=3.55 $ and $ x_{0}=0.6 $, which orbit may converge to fixed point of period 8.
\end{itemize}

We used Matlab with double precision to compute the calculations. The first interval of each case is obtained rounding $ x_{0} $ up and down. In other words, its endpoints $ \underline{x_{0}} $ and $ \overline{x_{0}} $ are the nearest floating points of $ x_{0} $.

\section{RESULTS}



The results obtained for each case with the calculations of $ F_{1} $ and the intersection $ F_{1}\cap F_{2} $ are presented in Tables \ref{tab1}, \ref{tab2} and \ref{tab3}. The fixed point, the width of the interval and the number of iterations reached by each method can be analysed at these tables.

\begin{table}[h]
	\caption{Solution of case 1 for $F_{1} $ and $ F_{1}\cap F_{2} $.}
	\label{tab1}
	\centering
	\begin{tabular}{l c c}	
		\hline
							 & $ F_{1} $ & $ F_{1}\cap F_{2} $ \\ \hline 
		Period of Fixed Point&$2$&$2$\\
		Width of Interval    &$3.65900\times10^{-6}$&$3.65833\times10^{-6}$\\
		Number of Iterations &$20$&$20$\\ \hline	
	\end{tabular}
\end{table}

\begin{table}[h]
	\caption{Solution of case 2 for $F_{1} $ and $ F_{1}\cap F_{2} $.}
	\label{tab2}
	\centering
	\begin{tabular}{l c c}	
		\hline
		& $ F_{1} $ & $ F_{1}\cap F_{2} $ \\ \hline 
		Period of Fixed Point&$4$&$4$\\
		Width of Interval    &$6.44715\times10^{-4}$&$4.72274\times10^{-4}$\\
		Number of Iterations &$23$&$23$\\ \hline	
	\end{tabular}
\end{table}

\begin{table}[h]
	\caption{Solution of case 3 for $F_{1} $ and $ F_{1}\cap F_{2} $.}
	\label{tab3}
	\centering
	\begin{tabular}{l c c}	
		\hline
							 & $ F_{1} $ & $ F_{1}\cap F_{2} $\\ \hline 
		Period of Fixed Point&$8$&$8$\\
		Width of Interval    &$9.14991\times10^{-4}$&$7.91852\times10^{-4}$\\
		Number of Iterations &$23$&$23$\\ \hline	
	\end{tabular}
\end{table}

In Tables \ref{tab1}, \ref{tab2} and \ref{tab3} we note that, despite the same number of iterations to converge to the same fixed point, there is a difference between the final sizes of intervals obtained by each method in each case.

\section{DISCUSSION}

By analysing the results, we note that the two calculations converge to the results expected according to \cite{Fei1978}. It shows that both were effective for these cases. The difference between the size of intervals obtained by these methods tends to grow when we approximate of unstable zones of (\ref{log}) with $ r $ approaching 4. Table \ref{wid} shows the sizes of intervals obtained for the first 10 iterations for each method to equation (\ref{log}) with $ r=3.6 $ and $ x_{0}=0.6 $. We can see that $ F_{1}\cap F_{2} $ has a lower growth rate compared to $ F_{1} $.

\begin{table}[h]
	\caption{Size of intervals obtained by $ F_{1} $ and $ F_{1}\cap F_{2} $.}
	\label{wid}
	\centering
	\begin{tabular}{c c c c}	
		\hline
	Iteration & $ F_{1} $            & $ F_{1}\cap F_{2} $  & Difference\\ \hline 
	$ 0 $     &$ 2.22\times10^{-16} $&$ 2.22\times10^{-16} $&$ 0 $\\
	$ 1 $     &$ 6.66\times10^{-16} $&$ 6.66\times10^{-16} $&$ 0 $\\ 
	$ 2 $     &$ 2.55\times10^{-15} $&$ 2.11\times10^{-15} $&$ 4.44\times10^{-16} $\\ 
	$ 3 $     &$ 9.77\times10^{-15} $&$ 7.77\times10^{-15} $&$ 2.00\times10^{-15} $\\
	$ 4 $     &$ 3.48\times10^{-14} $&$ 2.77\times10^{-14} $&$ 7.10\times10^{-15} $\\
	$ 5 $     &$ 1.25\times10^{-13} $&$ 9.99\times10^{-14} $&$ 2.60\times10^{-14} $\\
	$ 6 $     &$ 4.53\times10^{-13} $&$ 3.60\times10^{-13} $&$ 9.36\times10^{-14} $\\
	$ 7 $     &$ 1.63\times10^{-12} $&$ 1.29\times10^{-12} $&$ 3.37\times10^{-13} $\\
	$ 8 $     &$ 5.87\times10^{-12} $&$ 4.66\times10^{-12} $&$ 1.21\times10^{-12} $\\
	$ 9 $     &$ 2.11\times10^{-11} $&$ 1.67\times10^{-11} $&$ 4.37\times10^{-12} $\\
	$ 10 $    &$ 7.62\times10^{-11} $&$ 6.04\times10^{-11} $&$ 1.57\times10^{-11} $\\ \hline
	\end{tabular}
\end{table}

\section{CONCLUSIONS}

Both methods are effective to estimate the values of fixed points. Our proposed method gets some slight better results as the size of intervals obtained by the method based on the intersections of the interval extensions is smaller or equal than that obtained by traditional approach. It means that the proposed method can be used to increase the accuracy of the solutions. For example, in Table 2, there is a reduction of approximately 26\% and in Table 3, this reductions is 13\%.

We expect that using more interval extensions the results may be improved, which is as task for future investigations. We also intend to analyse the behaviour of the method in chaotic regime and examine how its use is feasible to help in the bifurcation diagram building for the logistic map and other maps. 

\section*{ACKNOWLEDGEMENTS}

We would like to thank CAPES, CNPq/INERGE, FAPEMIG and Federal University of S\~ao Jo\~ao del-Rei by support.

\bibliographystyle{ieeetr}
\bibliography{refer_nsc}

\begin{thebibliography}{10}

\bibitem{Galias2013}
Z.~Galias, ``The dangers of rounding errors for simulations and analysis of
  nonlinear circuits and systems? and how to avoid them,'' {\em Circuits and
  Systems Magazine, IEEE}, vol.~13, no.~3, pp.~35--52, 2013.

\bibitem{Nep2014}
E.~G. Nepomuceno, ``Convergence of recursive functions on computers,'' {\em The
  Journal of Engineering}, vol.~1, no.~1, 2014.

\bibitem{Nepomuceno2016}
E.~G. Nepomuceno and S.~A.~M. Martins, ``A lower-bound error for free-run
  simulation of the polynomial narmax,'' {\em Systems Science \& Control
  Engineering}, no.~ja, pp.~1--14, 2016.
\newblock Doi: 10.1080/21642583.2016.1163296.

\bibitem{HYG1987}
S.~Hammel, J.~Yorke, and C.~Grebogi, ``Do numerical orbits of chaotic dynamical
  processes represent true orbits?,'' {\em Journal of Complexity}, vol.~3,
  no.~2, pp.~136--145, 1987.

\bibitem{May1976}
R.~M. May, ``Simple mathematical models with very complicated dynamics,'' {\em
  Nature}, vol.~261, pp.~459--467, 1976.

\bibitem{MM1979}
R.~E. Moore, {\em Methods and applications of interval analysis}, vol.~2.
\newblock SIAM, 1979.

\bibitem{Ruetsch2005}
G.~Ruetsch, ``An interval algorithm for multi-objective optimization,'' {\em
  Structural and Multidisciplinary Optimization}, vol.~30, no.~1, pp.~27--37,
  2005.

\bibitem{alefeld2012}
G.~Alefeld and J.~Herzberger, {\em Introduction to interval computation}.
\newblock Academic press, 2012.

\bibitem{Rump1999}
S.~M. Rump, {\em INTLAB Interval Laboratory}.
\newblock Springer, 1999.

\bibitem{Rodrigues-Junior2015}
H.~M. Rodrigues~Junior and E.~G. Nepomuceno, ``Uso da computa\c{c}\~{a}o por
  intervalos para c\'{a}lculo de ponto fixo de um mapa discreto.,'' in {\em
  Proceedings of DINCON 2015 - Confer\^{e}ncia Brasileira de Din\^{a}mica,
  Controle e Aplica\c{c}\~{o}es}, 2015.

\bibitem{Fei1978}
M.~J. Feigenbaum, ``{Quantitative Universality for a Class of Non-linear
  Transformations},'' {\em {journal Of Statistical Physics}}, vol.~{19},
  no.~{1}, pp.~25--52, {1978}.

\end{thebibliography}

\end{document}